\documentclass[3p,times]{article}
\usepackage{amsfonts}
\usepackage{amsmath}
\usepackage[all,cmtip]{xy}
\usepackage{amssymb}
\usepackage{tikz-cd}
\usepackage{amscd}
\usepackage{cancel}
\usepackage{tabularx}
\usepackage{mathtools}
\usepackage{here}
\usepackage{lipsum}
\usepackage{lineno}
\usepackage{tikz}
\newenvironment{proof}{\noindent{\bf Proof.}}
{\noindent \ \hfill$\Box$\par}
\usepackage[all]{xy}
\usepackage{xcolor}
\newcommand{\bzd}{\hfill\boxed{}}
\usepackage{bbm}
\usepackage{colortbl}
\usepackage{selinput}
\usepackage{hyperref}
\hypersetup{
    colorlinks=true,
    linkcolor=blue,
    filecolor=magenta,
    urlcolor=cyan,
    pdftitle={Overleaf Example},
    pdfpagemode=FullScreen,
    }
\usepackage[a4paper, left=1in, right=1in, top=1in, bottom=1in]{geometry}
\newtheorem{theorem}{Theorem}[section]
\newtheorem{definition}[theorem]{Definition}
\newtheorem{pro}[theorem]{Proposition}
\newtheorem{lem}[theorem]{Lemma}
\newtheorem{cor}[theorem]{Corollary}
\newtheorem{rem}[theorem]{Remark}
\newtheorem{exam}[theorem]{Lemma}
\newcommand{\hc}{\circ}
\newcommand{\n}{\{ }
\newcommand{\nn}{\} }
\newcommand{\ca}{\eta}
\newcommand{\cg}{\sigma}
\newcommand{\lt}{\varepsilon}
\newcommand{\ch}{ \bar{\nu}  }

\newcommand{\ck}{\zeta}
\newcommand{\cn}{\kappa}
\newcommand{\co}{\rho}

\newcommand{\cs}{\bar{\zeta}}
\newcommand{\ccs}{\bar{\sigma}}
\newcommand{\ct}{\bar{\kappa}}

\newcommand{\w}{\omega}
\newcommand{\wq}{\infty}
\newcommand{\af}{\alpha}

\newcommand{\Ker}{\mathrm{Ker}}
\newcommand{\cok}{\mathrm{Coker}}

\newcommand{\ind}{\mathrm{Ind}}

\newcommand{\m}{\;\mathrm{mod}\,}
\newcommand{\pa}{\partial}

\newcommand{\e}{\iota}

\newcommand{\dyd}{\supseteq}
\newcommand{\xyd}{\subseteq}

\newcommand{\ty}{ \equiv}

\newcommand{\lo}{\Omega}
\newcommand{\tg}{\approx}

\newcommand{\jia}{\oplus}
\newcommand{\z}{\mathbb{Z}}

\newcommand{\s}{\Sigma}
\newcommand{\is}{\cong}

\newcommand{\mtodaa}[5]{\left\{
#1,\hspace{0.5mm}%
\begin{array}{c}
#2 \\
#4
\end{array}\hspace{-1.5mm},\hspace{0.5mm}%
\begin{array}{c}
#3 \\
#5
\end{array}\hspace{-1.5mm}
\right\}}
\usepackage{array}

\newcolumntype{Y}{>{\centering\arraybackslash}X}

\newcommand{\q}{ \Delta}

\newcommand{\mtoda}[5]{\left\{
\begin{array}{c}%
{#1} \\
{#4}
\end{array}\hspace{-1.5mm},\hspace{0.5mm}%
\begin{array}{c}%
{#2} \\
{#5}
\end{array}\hspace{-1.5mm},\hspace{0.5mm}%
{#3}
\right\}}
\begin{document}

\date{}
\title{An Unstable Approach to   the May-Lawrence   Matrix Toda bracket and the \textit{2}nd James-Hopf Invariant}
\author{Juxin Yang\thanks{\text{Beijing Institute of Mathematical Sciences and Applications, Beijing, 101408, P.R. China.  yangjuxin@bimsa.com}},  \,\;Toshiyuki Miyauchi\thanks{Department of Applied Mathematics, Faculty of Science, Fukuoka University;
8-19-1 Nanakuma, Jonan-ku, Fukuoka 814-0180, Japan}\;  and  Juno Mukai \thanks{Shinshu University;
3-1-1 Asahi, Matsumoto, Nagano 390-8621, Japan;
jmukai@shinshu-u.ac.jp}}
\maketitle
\begin{abstract}
In this paper, we  give an unstable approach of   the May-Lawrence   matrix Toda bracket, which becomes a useful tool for the theory of  determinations of unstable homotopy groups. Then, we  give a generalization of the classical isomorphisms between  homotopy groups of $(JS^{m},S^{m}) $ and $(JS^{2m},*)$  localized at 2. After that we  provide a  generalized    $H$-formula for matrix Toda brackets. As an application, we show a new construction of $\ct'\in\pi_{26}(S^{6})$ localized at 2 which improves the construction of  $\ct'$  given by \cite{20STEM}.
\end{abstract}

\section{Introduction}

The Toda bracket serves the art on constructing homotopy liftings and extensions of maps, it plays a fundamental role on dealing with  composition relations of homotopy classes. The theory of matrix  Toda brackets  is widely studied, there are both classical and category-theoretic descriptions for such  Toda brackets. M.G. Barratt (\cite[1963]{Bar}) first studies the   generalized  Toda bracket which is  called the matrix Toda bracket nowadays, he uses it to   detect some elements of homotopy groups of spheres. Then, M. Mimura (\cite[1964]{M}) introduces Barratt's matrix Toda bracket and gives some important  properties of it.  K.A. Hardie, H.J. Marcum and N. Oda (\cite[2001]{Oda2001}) define the box bracket which greatly generalized the 3-fold Toda bracket.  They also give applications of the box bracket for unstable homotopy groups of spheres. Recently,  the second, third authors\,(\cite[2017]{32STEM}) define the left matrix Toda brackets  indexed by $n$ for all $n\geq0$. They use such tools to define an element  $\bar{\ccs}_{14}\in\pi_{48}(S^{14};2)$ with nontrivial Hopf invariant, which is one of the key ingredients for the project  on  determinations of $\pi_{32+k}(S^{k};2)$ for all $k\geq2$\;(\cite{32STEM}).

 Matrix Toda brackets of the most general form and of the largest shape are proposed by Peter May (\cite{May}) and established by his student A.F. Lawrence (\cite{Mays})  in 1969.
 Such matrix Toda brackets are of great importance in the stable homotopy theory, see \cite[p. 158]{wgzxzl}. However, they do not define the   matrix Toda brackets indexed by \textit{n}, ($n\geq1$). In the unstable range, their  matrix Toda brackets  are sometimes too large, resulting in an inadequate exploitation of the available information.   As readers will observe in Remark \ref{xbh}, the presence of the indexing number $n\geq1$ in a Toda bracket is crucial in the unstable homotopy theory. Hence, in this paper, 
we  give an unstable approach  to   the May-Lawrence   Matrix Toda brackets  study the properties. 

After localization at 2, by the methods of  relative homotopy theory,  we generalize the classical group isomorphisms 
 $\pi_{r}(JS^{m}, S^{m})\stackrel{H_{2*}\,}\longrightarrow \pi_{r}(JS^{2m},*)$ to  isomorphisms $[ (C\s^{r-1}Z, \s^{r-1}Z), (JS^{m}, S^{m})]\stackrel{H_{2*}\;\;}\longrightarrow [(C\s^{r-1}Z, \s^{r-1}Z), (JS^{2m},*)]$ for any pointed  \textit{CW} complex $Z$. Successively, we provide a generalization of Toda’s formula in \cite[Proposition\,2.6]{Toda} for the
value of a matrix Toda bracket under the homomorphism $H$ of the \textit{EHP} sequence. Such
 a formula for left matrix Toda brackets is missing in \cite[Section 4]{32STEM}. As an application, we give a new construction of $\ct'\in\pi_{26}(S^{6})$ localized at 2  for \cite{20STEM}, by which we are able to derive $H\n\nu_{6}, \;\ca_{9}, \;\ca_{10}\cn_{11} \nn= \ca_{11}\cn_{12}$.  Attempts to obtain  such a formula  by classical Toda bracket methods have not succeeded. By this construction we can further obtain  $\q\ct'=[\nu_{5}]\nu_{8}\cn_{11}\in\pi_{25}(\textit{SO}(6))_{(2)}$, 
which is an important ingredient in the determination of $\pi_{25}(\textit{SO}(n))_{(2)}$ 
for $n\geq3$; we hope to provide the details of this calculation in a future paper.\\\indent 
Herein lies our first main theorem.

\begin{theorem}\label{hao1} 

After localization at 2, suppose $Z$ is a pointed  \textit{CW} complex and $H_{2}: JS^{m}\rightarrow JS^{2m}$
      is the second James-Hopf invariant. Then there exist   isomorphisms of $\z_{(2)}$-modules for all $r\geq2$,  $$H_{2*}:[(C\s^{r-1}Z, \s^{r-1}Z),\;(JS^{m},S^{m})]\stackrel{\is}\longrightarrow[(C\s^{r-1}Z, \s^{r-1}Z), (JS^{2m},*)].$$  
\end{theorem}

We recall the classical $H$-formula given by Toda. 
Let spaces be localized at 2. For the suspended homotopy classes $\s\af,\s\beta$ and $\s\gamma$ with respect to the sequence\vspace{-0.67\baselineskip}
$$
\xymatrix@C=1cm{
 S^{m+1}\ar@{<-}[r]^{\s \af}&\s K \ar@{<-}[r]^{\;\,\s\beta} & S^{k+1}\ar@{<-}[r]^{\s \gamma}& S^{i+1},} $$
in \cite[Proposition\,2.6]{Toda}, Toda gives a formula which explicates the interaction of the homomorphism $H$ with $P$ in the process of using the  Toda bracket and the \textit{EHP} sequence. That is, $$H\n \s\af,\s\beta,\s\gamma\nn_{1}=-P^{-1}(\af\hc\beta)\hc\s^{2}\gamma,$$ which provides one of the most  fundamental strategies to construct new elements of homotopy groups of spheres. In this article, we generalize it to the matrix Toda bracket. It is noteworthy that  our $H$-formula for the matrix Toda bracket  cannot be straightforwardly  derived from Toda's $H$-formula, but rather deduced through the  theory of  groups $[(C\s X,\s X), (Y, B)]$.

Our second main theorem   is stated as follows.

\begin{theorem}\label{zdl2} 

Localized at 2, let $m,\ell,r\in\z_{+}$  and the homotopy classes
$f_{k}$, $b_{k,s}$ together with $c_{s}$ ($1\leq k\leq \ell , 1\leq s\leq r$)  be given with respect to the following sequences,\\\\\centerline{
$
\xymatrix@C=1cm{
 S^{m+1}\ar@{<-}[r]^{\s f_{k}}&\s X_{k} \ar@{<-}[r]^{\;\,\s^{}b_{k,s}} & \s^{}Y_{s}\ar@{<-}[r]^{\s^{} c_{s}}& \s^{}U,} $}
\;\\ 
\noindent where all homotopy classes $c_{s}$ are suspensions.
Suppose the formal matrices  $$\mathbbm{f}_{1\times\ell}=(f_{1},f_{2},\cdots,f_{\ell}),\;\; \mathbbm{b}_{\ell\times r}=(b_{k,s}) \;\;\text{and}\;\; \mathbbm{c}_{r\times 1}=(c_{1},c_{2},\cdots,c_{r})^{T}$$ satisfy $\s(\mathbbm{f}_{1\times\ell}\hc\mathbbm{b}_{\ell\times r})= O $ and $\mathbbm{b}_{\ell\times r}\hc\mathbbm{c}_{r\times 1}= O$. Write $\mathbbm{b}_{\bullet, s}=(b_{1,s},b_{2,s}\cdots, b_{\ell,s})^{T}$, the  $s$-$th$ column   vector of $\mathbbm{b}_{\ell\times r}$.
Then,$$H\{\s \mathbbm{f}_{1\times\ell}, \s\mathbbm{b}_{\ell\times r}, \s\mathbbm{c}_{r\times 1}\}_{1}=-\sum_{s=1}^{r}P^{-1}(\mathbbm{f}_{1\times\ell}\hc\mathbbm{b}_{\bullet,s})\hc\s^{2}\mathbbm{c}_{r\times 1}. $$

\end{theorem}

\noindent\textbf{Acknowledgement} The authors  are indebted  to  professor Dai Tamaki,  Sergei O. Ivanov  and Yoshihiro Hirato for   many fruitful conversations on
this project.\,The first author would like to thank professor Dai Tamaki very much for invitation to Shinshu University by his grant JSPS KAKENHI\,(No.\,20K03579). And
the second author is supported in part by JSPS KAKENHI\,(No.\,22K03326). 
\section{Preliminaries}

\subsection{Notations}
 In this paper, all spaces, maps, and homotopy classes are pointed with $*$ denoting  base points unless otherwise stated. Thus $(X, A)$ is an abbreviation for the triad
$(X, A, *)$, and so on. By  a ``space" we mean a  \textit{CW} complex,   a sphere  $S^{m}$ implicitly has $m\geq1$.  When saying a homotopy class $\af$ being a suspension, namely, $\af=\s\af'$ for some homotopy class $\af'$, it is invariably assumed that both the domain and codomain of $\af$ are  suspensions.

 Whenever it does not cause misunderstanding, in this article the phrase  ``unstable''  will refer to ``not necessarily stable",  following the convention used in the unstable Adams spectral sequence.

The symbol $C(-)$  stands for the reduced cone functor. For a map $f:X\rightarrow Y$,   use $C^{f}:CX\rightarrow CY$  to denote the extended map over $f$ by applying  the reduced cone functor. Essentially,   $C^{f}$ is $f\wedge \mathrm{id}_{I}$. As usual, $C_{f}$  denotes the homotopy cofibre of $f$. 

We  identify a set of a single element with its element if no confusions arise. 

The notation $\z_{(2)}$  denotes the group or the ring of  2-local integers. We use  $P$ to denote the symbol $\Delta$ in \cite{Toda},  the boundary homomorphism of the \textit{EHP} sequence. 

Suppose $U$ is an abelian group and $V$ is a subgroup. If we have a coset $S\in U/V$, then $V$ is often called the  \textit{\textbf{indeterminacy}} of $S$, denoted by $\mathrm{Ind}(S)=V$.

\subsection{Some fundamental facts}\label{name}
Given the following sequence of spaces and homotopy classes with   $Z$  a suspension such that $\af\hc \s^{n}\beta=\beta\gamma=0$,\vspace{-0.7\baselineskip}\;\\\;
\[\xymatrix@!C=0.83cm{ &W&  \s^{n}X\; & \s^{n}Y  &\s^{n}Z,& 
       \ar"1,4";"1,3" _{\;\s^{n}\beta_{}} \ar"1,5";"1,4" _{ 
 \s^{n}\gamma}
\ar"1,3";"1,2" _{\af}
     }\]
the (3-fold) Toda bracket  which is denoted by  $\n\af,\s^{n}\beta,\s^{n}\gamma\nn_{n}$ is defined to be the set of all compositions of the form \vspace{-0.3\baselineskip}$$(-1)^{n}\text{ext}_{_{\s^{n}\beta}}(\af)\hc \s^{n}\text{coext}_{_{\beta}}(\gamma),$$ where  $\text{ext}_{_{\s^{n}\beta}}(\af)\in [C_{\s^{n}\beta},\,W]$ is an extension of $\af$ with respect to  $\s^{n}\beta$, and $\text{coext}_{_{\beta}}(\gamma)\in[\s Z, C_{\beta}] $ is a coextsion of $\gamma$ with respect to  $\beta$, (in \cite[p. 13]{Toda}, $\text{ext}_{_{\s^{n}\beta}}(\af)$ and $\text{coext}_{_{\beta}}(\gamma)$ are denoted by $\overline{\af}$ and $\widetilde{\gamma}$  respectively; for more properties of extensions and coextensions, see \cite[Lemma   2.2 and Lemma   2.3]{dd}). Such a Toda bracket is a coset included in $[\s^{n+1}Z,W]$, that is,  $$\n\af,\s^{n}\beta,\s^{n}\gamma\nn_{n}\in [\s^{n+1}Z,W]/A,$$ where \vspace{-1.85\baselineskip} \\\;\begin{equation}\label{22a}
    A=\af\hc \s^{n}[\s Z,X]+[\s^{n+1}Y, W]\hc \s^{n+1}\gamma. \tag{2.2a}
\end{equation}

\noindent Furthermore, as a group homorphism sends a coset to a coset, it is clear that $F(\n\af,\s^{n}\beta,\s^{n}\gamma\nn_{n})$
is a coset of $F(A)$ where $F$ is a group homomorphism; in particular,
since the composition operators $f\hc(-)=f_{*}$ and $(-)\hc (\s g)=(\s g)^{*}$ give homomorphisms, we derive
$f\hc \n\af,\s^{n}\beta,\s^{n}\gamma\nn_{n}$
is a coset of $f\hc A$,\;and $\n\af,\s^{n}\beta,\s^{n}\gamma\nn_{n}\hc \s g$\;
      is a coset of $A\hc \s g$. For a fixed $n$, the above Toda bracket is a function of the triad $(\af, \beta,\gamma)$; more specifically speaking, for a given  integer $n$, the  Toda bracket $\n \af,\s^{n}\beta,\s^{n}\gamma\nn_{n}$   depends only on $(\af, \beta,\gamma)$ but not $(\af, \s^{n}\beta,\s^{n}\gamma)$.  It is necessary to point out that  even if\; $\s^{n}\beta'=\s^{n}\beta, \s^{n}\gamma'=\s^{n}\gamma \; \text{and}\; \beta'\gamma'=0$,
      $$\n \af,\s^{n}\beta',\s^{n}\gamma'\nn_{n}\neq\n \af,\s^{n}\beta,\s^{n}\gamma\nn_{n}\;\text{in general},$$
  (see \cite[Remark 3.1]{dd}).
   For more basic properties of Toda brackets, see \cite[p. 10-12]{Toda}.

   It is noteworthy that the  Toda bracket indexed by $n$ with $n\geq1$ makes better use of the desuspension properties of homotopy classes and plays a critical role in the theory of computing unstable homotopy groups. By presenting the following three arguments, we substantiate  in the unstable range the essential significance of the condition the indexing number $n\geq1$.\begin{rem}\label{xbh}
     
\begin{itemize} 
 \item [(i)]   It proves to be advantageous for defining new elements in $\pi_{*}(S^{i})$. For example, recall from Toda's work (\cite[Lemma 6.1, p. 51]{Toda}) that  $\lt_{3}\in\n\ca_{3},\s\nu', \nu_{7} \nn_{1}$. The presence of the subscript 1 enables us to utilize the $H$-formula, $H(\lt_{3})\in\n\ca_{3},\s\nu', \nu_{7} \nn_{1}=-P^{-1}(\ca_{2}\hc\nu')\hc\nu_{8}=\nu_{5}^{2}$, thereby  confirming $\lt_{3}$ to be a lifting of $\nu_{5}^{2}.$  This, in fact, constitutes the fundamental method for constructing new elements of homotopy groups of spheres in \cite{Toda}. The vast majority of new elements in \cite{Toda} are constructed using this method.

 \item [(ii)]   It is useful to consider the connecting homomorphism. Suppose we aim to determine $\pi_{k}(X)$, with $\pi_{*}(Y)$ and $\pi_{*}(F)$ being already known where $F\rightarrow X \rightarrow Y$ is a homotopy fibration.
    We need consider the exact sequence
\begin{tikzcd}\
&&&\pi_{k+1}(Y) \arrow[r,"\pa_{_{(k+1)}}"] & \pi_{k}(F) \arrow[r] & \pi_{k}(X)\arrow[r] & \pi_{k}(Y) \arrow[r,r,"\pa_{_{(k)}}"] & \pi_{k-1}(Y) 
\end{tikzcd}

inducing a short exact sequence\\ \begin{tikzcd}
&&&0\arrow[r,] & \cok(\pa_{_{(k+1)}}) \arrow[r] & \pi_{k}(X)\arrow[r] & \mathrm{Ker}(\pa_{_{(k)}}) \arrow[r] & 0.
\end{tikzcd}
 
  \noindent Then for $\mathbbm{x}\in\n\af, \s^{n}\beta,\s^{n}\gamma \nn_{n}\xyd\pi_{k+1}(Y)$, the relation $$\pa_{_{(k+1)}}(\mathbbm{x})\in \pm \n \pa_{_{(k+1)}}(\af), \s^{n-1}\beta,\s^{n-1}\gamma \nn_{n-1}$$   holds for $n\geq1$ but  not for $n=0$. Similarly, considering the image of $\mathbbm{x}\in\n\af, \s^{n}\beta,\s^{n}\gamma \nn_{n}\xyd\pi_{k+2}(S^{2m+1}) $ under the homomorphism $P:\pi_{k+2}(S^{2m+1})\rightarrow\pi_{k}(S^{m}) $    in the 2-local   EHP sequence, the   relation $$P(\mathbbm{x})\in \pm \n P(\af), \s^{n-2}\beta,\s^{n-2}\gamma \nn_{n-2}$$   holds for $n\geq2$ but  not for $n\in\n0, 1\nn$.

\item [(iii)]   During studying unstable homotopy groups, it can  trace back to the dimension before the birth of certain troublesome elements, thereby reducing interference. For instance, notice the proof of determining $P(\lt_{13})$, (\cite[p. 78]{Toda}). It tells us that localized at 2, $P(\lt_{13})\in \n \nu_{6},\ca_{9}, 2\e_{10}\nn_{6}\hc(-\lt_{11})=\nu_{6}\hc \s^{6}\n\ca_{3},2\e_{4}, \lt_{4}
 \nn\xyd \nu_{6}\hc\s^{6}\pi_{13} (S^{3})$. Computing $\nu_{6}\hc\s^{6}\pi_{13} (S^{3})$ is easier than 
computing $\nu_{6}\hc\pi_{19} (S^{9})$ due to the fact that some unstable elements in $\pi_{13} (S^{3})$ are killed by the iterated suspension functor $\s^{6}$. An alternative perspective is that, Formula (\ref{22a}) reveals it is typically more practical to calculate the indeterminacy of the Toda bracket indexed by $n\geq1$. As   some troublesome elements usually appear in $[\s^{n+1} Z,\s^{n} X]$ but not in $[\s Z, X]$.

\end{itemize}
  \end{rem}

   \indent We shall  give a brief introduction of Toda's naming convention  for the  generators of $\pi_{*}(S^{n};\,2\,)$ to help readers to read some of our propositions and their proofs conveniently. This naming convention is also adopted by \cite{20STEM}, \cite{2122STEM},\cite{2324STEM}, \cite{Oda}, \cite{32STEM} and so on. Since Toda's 2-primary component method  in \cite{Toda}  naturally corresponds with the 2-localization method by which we can deal with spaces more intuitively, in this article  we prefer to use the language  of the 2-localization to state the results on  the 2-primary components.\\\indent Let spaces be localized at 2. Roughly speaking, 
 for $\mathbbm{x}$ 
 which represents a Greek letter,   $\mathbbm{x}_{n}$   denotes one of the  generators of $\pi_{n+r}(S^{n})$ for some $r$, the  subscript $n$ indicates the codomain of $\mathbbm{x}_{n}$. Moreover,  $\mathbbm{x}_{n+k}:=\s^{k}\mathbbm{x}_{n}$, $\mathbbm{x}:=\s^{\wq}\mathbbm{x}_{n}$, and $\mathbbm{x}_{n}^{\ell}$ is the abbreviation of $\mathbbm{x}_{n}\hc \mathbbm{x}_{n+r}\hc \cdots \hc \mathbbm{x}_{n+(\ell-1)r}$,\,\, ($\ell$ factors). The usages of  $\overline{\mathbbm{x}}_{n} $ and $\mathbbm{x}^{*}_{n}$ are similar to  above. In $\pi_{j+r}(S^{j})$ (not a stable homotopy group), if a  generator is written  without a subscript,  then  this  generator  does not survive in the stable homotopy group $\pi_{r}^{S}(S^{0})$ or its  $\s^{\wq}$-\,image is divisible by 2. For instance, for  $\theta\in\pi_{24}(S^{12})$,\;$\cg'''\in\pi_{12}(S^{5})$, their  $\s^{\wq}$-\,images satisfy $\s^{\wq}\theta=0$,   $\s^{\wq}\cg'''=8\cg$   of order 2. There is an  advantage of this naming convention, that is, we can examine the commutativity of unstable compositions conveniently,  $$\mathbbm{x}_{n}\hc \mathbbm{y}_{i}=\pm \mathbbm{y}_{n}\hc\mathbbm{x}_{j}  \;\;\text{ for some}\; i, j\;\;\; \text{if} \;\;n\geq a+b,$$
 where   $\n \mathbbm{x}_{k}\nn$ was born in $\pi_{*}(S^{a})$ and   $\n \mathbbm{y}_{\ell}\nn$ was born in $\pi_{*}(S^{b})$, (see \cite[Proposition 3.1]{Toda}). For example, for the  elements  in \cite{Toda},  $$\cg_{n}\in\pi_{7+n}(S^{n})\,(n\geq8) \;\text{and} \;\mu_{n}\in\pi_{9+n}(S^{n})\,(n\geq3),$$ we have  $\cg_{8+3}\mu_{i}=\pm \mu_{8+3}\cg_{j}$, successively,\; $\cg_{11}\mu_{18}=\mu_{11}\cg_{20}$, ($\pm$ is not necessary, since $\mu_{3}$ is of order 2). But  $\cg_{10}\mu_{17}\neq\mu_{10}\cg_{19}$, (see \cite[p. 156]{Toda}).
 Some common  generators are summarized in \cite[p. 189]{Toda} and \cite[(1.1),\,p. 66]{Oguchi}.

\section{An unstable approach to   the May-Lawrence   Matrix Toda bracket}\label{333}

We will give an unstable approach to   the May-Lawrence   Matrix Toda bracket. Recall our Remark \ref{xbh};  in the unstable range, the May-Lawrence   matrix Toda brackets are sometimes too large, resulting in an inadequate exploitation of the available information.   We need to reconsider them in the unstable range and study them  indexed by $n$ for all  $n\geq0$, in particular, $n\geq1$.

This reconsideration is  necessary. For example, in Proposition \ref{ind}, we point out that the condition $n\geq1$ is required for the formula to calculate the indeterminacy. More importantly, we  introduce the generalized $H$-Formula, (Theorem \ref{zdl2}), a crucial formula for defining new elements in the unstable homotopy  group $\pi_{*}(S^{i}_{(2)})$. This generalized $H$-Formula cannot be directly derived from the $H$-Formula in \cite[Proposition 2.6]{Toda}. The establishment necessitates a comprehensive grasp of the intricate proof of  the classical $H$-Formula, alongside the application of the  theory of  groups $[(C\s X,\s X), (Y, B)]$.

Suppose $f_{i}: T_{i}\rightarrow W\; (i=1,2)$ are maps or homotopy classes. By abuse of notation, the following composition is still denoted by $f_{1}\vee f_{2}$,\\
\xymatrix@C=3em{
&&&& W \ar@{<-}[r]^{\nabla
 \;\;}_{\text{folding}\;\;} & W\vee  W \ar@{<-}[r]^{\; f_{1}\vee f_{2}} & T_{1}\vee  T_{2}. 
}

\noindent And we use the identification $\s^{n}(V_{1}\vee V_{2})=\s^{n}V_{1}\vee \s^{n}V_{2},\;\,  (V_{1}, V_{2}\colon$ spaces). For  $k$-fold wedges, we make comparable assumptions.

Given $\ell,r\in\z_{+}$, let $W, X_{k}, Y_{s}$ and $U$ be spaces, ($k,s\in\z$,  $1\leq k\leq \ell,\;1\leq s\leq r$). For each $k$ and $s$, let the  homotopy classes
$a_{k}$, $b_{k,s}$ and $c_{s}$ be given with respect to the following sequences\\\\\centerline{
$
\xymatrix@C=1cm{
 W\ar@{<-}[r]^{a_{k}}&\s^{n}X_{k} \ar@{<-}[r]^{\;\,\s^{n}b_{k,s}} & \s^{n}Y_{s}\ar@{<-}[r]^{\s^{n} c_{s}}& \s^{n}U.} $}
\;\\ 
\noindent Moreover, suppose $\bigvee_{p=1}^{\ell} X_{p}\stackrel{\mathrm{i}_{k}}\longleftarrow X_{k}$ and $\bigvee_{q=1}^{r} Y_{q}  \stackrel{\;\mathrm{j}_{s}}\longleftarrow Y_{s}$ are the inclusions.
We introduce formal matrices to denote homotopy classes,\vspace{-1.15\baselineskip}
$$\mathbbm{a}_{1\times \ell}=( 
a_{1}, a_{2}, \cdots, a_{\ell} ):=a_{1}\vee a_{2}\vee\cdots\vee a_{\ell}\in [\s^{n}(\bigvee_{p=1}^{\ell} X_{p}),\; W ],$$\vspace{-0.61\baselineskip}
$$\mathbbm{b}_{\bullet,s}=(b_{1,s},\; b_{2,s},\; \cdots, \;b_{\ell,s})^{T}:=\mathrm{i}_{1}\hc b_{1,s}+\mathrm{i}_{2}\hc b_{2,s}+\cdots+\mathrm{i}_{\ell}\hc b_{\ell,s}\in [Y_{s}, \bigvee_{p=1}^{\ell} X_{p}],$$
\vspace{-0.31\baselineskip}
$$\mathbbm{b}_{\ell\times r}=(b_{k,s}):=\mathbbm{b}_{\bullet,1}\vee \mathbbm{b}_{\bullet,2}\vee\cdots \vee \mathbbm{b}_{\bullet,r}\in [ 
    \bigvee_{q=1}^{r} Y_{q}   ,\bigvee_{p=1}^{\ell} X_{p}],$$\vspace{-0.31\baselineskip}
$$\s^{n}\mathbbm{b}_{\ell\times r}:=(\s^{n}b_{k,s})=\s^{n}(\mathbbm{b}_{\bullet,1}\vee \mathbbm{b}_{\bullet,2}\vee\cdots \vee \mathbbm{b}_{\bullet,r})\in [\s^{n}( 
    \bigvee_{q=1}^{r} Y_{q}  )  ,\;\s^{n}(\bigvee_{p=1}^{\ell} X_{p})]$$\vspace{-0.31\baselineskip}
$$\text{and}\qquad\mathbbm{c}_{r\times 1}=( c_{1}, c_{2}, \cdots, c_{\ell})^{T}:=\mathrm{j}_{1}\hc c_{1}+\mathrm{j}_{2}\hc c_{2}+\cdots+\mathrm{j}_{\ell}\hc c_{\ell}\in [U, \bigvee_{q=1}^{r} Y_{p}].\;\qquad\qquad\;$$
\noindent A homotopy class is regarded as a $1 \times 1$ formal matrix.   We say that two formal matrices are equal if the corresponding entries  are equal.

Through direct verification, it is evident that the composition of $\mathbbm{a}_{1\times \ell}$ and $\mathbbm{b}_{\ell\times r}$  as homotopy classes  corresponds to their multiplication as formal matrices; 
the composition of   $\mathbbm{b}_{\ell\times r}$ and $\mathbbm{c}_{r\times 1}$ as homotopy classes  corresponds to their multiplication as formal matrices if $c_{s}$ are suspensions for all $1\leq s\leq r$. (Notice that $f$ being a suspension is a sufficient condition for $f^{*}(-)=(-)\hc f $ being  a homomrphism). The formal matrices of trivial homotopy classes are denoted by $O$. For any formal matrix $\mathbbm{x}_{\ell\times r}$, the notation $\mathbbm{x}_{\bullet,s}$ always denotes its $s$-$th$ column  vector. 

Then, we are poised to present the ensuing definition.
 
\begin{definition}\label{lwndy}
Using above notations, further suppose the homotoy classes $c_{s}$ are suspensions for all $1\leq s\leq r$, $\mathbbm{a}_{1\times \ell}\hc\s^{n}\mathbbm{b}_{\ell\times r}=O  $
 and  $\mathbbm{b}_{\ell\times r}\hc\mathbbm{c}_{r\times 1}=O$. The matrix Toda bracket  indexed by $n$, $\n \mathbbm{a}_{1\times \ell},\; \s^{n}\mathbbm{b}_{\ell\times r}, \;\s^{n}\mathbbm{c}_{r\times 1}\nn_{n}$, is defined to be the \textit{3}-fold Toda bracket with respect to the following sequence\vspace{-0.2\baselineskip}$$\xymatrix@C=1.3cm{
 W\ar@{<-}[r]^{\mathbbm{a}_{1\times \ell}\qquad}&\s^{n}(\bigvee_{k=1}^{\ell} X_{k}) \ar@{<-}[r]^{\s^{n}\mathbbm{b}_{\ell\times r}} & \s^{n}( \bigvee_{s=1}^{\ell} Y_{s})\ar@{<-}[r]^{\qquad\s^{n}\mathbbm{c}_{r\times 1}}& \s^{n}U.} $$
    
\end{definition}
\noindent The matrix Toda bracket $\n \mathbbm{a}_{1\times \ell},\; \mathbbm{b}_{\ell\times r}, \mathbbm{c}_{r\times 1}\nn_{0}$ is denoted by $\n \mathbbm{a}_{1\times \ell},\; \mathbbm{b}_{\ell\times r}, \mathbbm{c}_{r\times 1}\nn_{}$ for short.

The following proposition presents an  efficient method for computing the indeterminacy of the matrix Toda bracket.

\begin{pro}\label{ind}
    Under the condition of Definition  \ref{lwndy},
 additionally  suppose $n\geq1$ and the spaces $X_{k}$ are suspensions for all $1\leq k\leq\ell$. Then, $\n \mathbbm{a}_{1\times \ell},\; \s^{n}\mathbbm{b}_{\ell\times r}, \;\s^{n}\mathbbm{c}_{r\times 1}\nn_{n}\in [\s^{n+1}U, W]/\mathrm{Ind},$ \;where\[\mathrm{Ind}=\sum_{k=1}^{\ell} (a_{k}\hc \s^{n}[\s U, X_{k}]) +       \sum_{s=1}^{r}    ([\s^{n+1}Y_{s}, X]\hc\s^{n+1}c_{s}).\]
 
\end{pro}
\begin{proof} Since there exist isomrphisms for any  path-connected spaces $A, B,C$, 
\begin{eqnarray}
& 
&\notag[\s A, \s B\vee \s C] \\\notag& 
\is& [A, \lo(\s B\vee \s C)]   \\\notag&
\is& [A, \lo \s B \times \lo \s C\times \lo \prod_{i=3}^{\wq} \s\Psi_{i}]  \\\notag&
\is& [\s A,  \s B] \times[\s A,  \s C] \times[\s A,  \prod_{i=3}^{\wq} \s\Psi_{i}]  , 
\end{eqnarray}

\noindent where $\prod_{i=3}^{\wq} \s\Psi_{i}$ is the infinite product   given by the \textit{Hilton-Milnor} theorem, (see \cite[Theorem  2.2]{jj}). The homomorphism $[\s A,  \prod_{i=3}^{\wq} \s\Psi_{i}]\hookrightarrow[\s A, \s B\vee \s C]$ is induced by iterated Whitehead products (\cite[Theorem  2.2]{jj},  \cite[Theorem  8.1,\,p. 533]{GW}) which are in $\Ker(\s)$. By assumption, we see that  $n\geq1$ and  the spaces $X_{k}$ are all  suspensions. Then, 
 we derive \begin{equation}
     (a_{1}\vee a_{2}\vee\cdots\vee a_{\ell})\hc\s^{n}[\s U,  \bigvee_{p=1}^{\ell} X_{p}]=\sum_{k=1}^{\ell} (a_{k}\hc \s^{n}[\s U, X_{k}]). \label{pj000}\tag{3.1a}
 \end{equation}  
\noindent
Therefore, the proposition  follows from Definition \ref{lwndy} and \cite[Lemma 1.1, p.  9]{Toda}.
\end{proof}
\;\\

For the above proposition, the condition  \,$n\geq1$\, is indispensable. Recall that we give some explanations on the names of the generators of $\pi_{*}(S^{n})$  in subsection \ref{name}.  We will freely use the notations in \cite{Toda} and \cite{20STEM}. 

The subsequent example demonstrates that for $n=0$, the formula given by  Proposition \ref{ind} does not hold in general.

\begin{exam}

After localization at 2, there is  a  matrix Toda bracket $$
R=         \n (\ca_{13},\cg_{13}), (\cg_{14},\ca_{20})^{T},4\ck_{21}\nn$$
\noindent with respect to the  sequence \vspace{-1.4\baselineskip}\;\\ 
\[\xymatrix@!C=1.7cm{&S^{13}&  S^{14}\vee S^{20}\;\; & \;\;S^{21}  &S^{32},& 
       \ar"1,4";"1,3" _{\,\,\quad\,\; \mathrm{j}_{1}\cg_{14}+\mathrm{j}_{2}\ca_{20}} \ar"1,5";"1,4" _{ 
 4\ck_{21}}
\ar"1,3";"1,2" _{\ca_{13}\vee \cg_{13}\quad}
     }\]\noindent where $\mathrm{j}_{1}: S^{14}\hookrightarrow S^{14}\vee S^{20}$  and $\mathrm{j}_{2}: S^{20}\hookrightarrow S^{14}\vee S^{20}$ are the inclusions. 
This matrix Toda bracket is a coset of\; 
$\mathrm{span}\n(\s\theta)\ch_{25},(\s\theta)\lt_{25}\nn\tg(\z/2)^{\jia2}$, but  not a coset of \;
$\ca_{13}\hc\pi_{33}(S^{14})+\cg_{13}\hc\pi_{33}(S^{20})=\mathrm{span}\n(\s\theta)\ch_{25}\nn\tg\z/2$.

\end{exam}
\begin{proof} Let spaces be localized at 2. Notice $4\cg_{10}\ck_{17}=0$,\,(\cite[(12.23),\,p. 163]{Toda}); then it is easy to check that $R$ is well-defined. By Definition \ref{lwndy}, we infer $R$ is a coset of $$A=(\ca_{13}\vee \cg_{13})\hc \pi_{33}(S^{14}\vee S^{20})+\pi_{22}(S^{13})\hc4\ck_{22}.$$ Since
$\pi_{22}(S^{13})\tg (\z/2)^{\jia3}$,  (\cite[p. 61]{Toda}), so  $\pi_{22}(S^{13})\hc4\ck_{22}=0$. Then  $A=(\ca_{13}\vee \cg_{13})\hc \pi_{33}(S^{14}\vee S^{20}).$

In order to compute the group $A$, we first undertake several preliminary calculations. By
\cite[Theorem  8.1,\,p. 533]{GW}, we have
$\pi_{33}(S^{14}\vee S^{20})=\mathrm{j}_{1*}(\pi_{33}(S^{14}))\jia \mathrm{j}_{2*}(\pi_{33}(S^{20}))\jia[\mathrm{j}_{1},\mathrm{j}_{2}]_{*}(\pi_{33}(S^{33})).$
 Through \cite[p. 162, p. 75]{Toda}, we see the homotopy groups are given by
$\pi_{33}(S^{14})=\mathrm{span}\n \w_{14}\nu_{30},\ccs_{14},\cs_{14}\nn\tg\z/2\jia\z/2\jia\z/8\;\text{and}\; \pi_{33}(S^{20})=0;$
following from   \cite[the proof of Lemma   (2.15),\,p. 75]{Oguchi} and \cite[Proposition\,3.1\,(3),\,p. 55]{Oda}, we obtain $\ca_{13}\w_{14}\nu_{30}=\lt_{13}^{*}\nu_{30}=(\s\theta)\ch_{25}.$
According to \cite[(2.7), p.  8]{2324STEM}, we infer $\ca_{6}\ccs_{7}=\ca_{6}\cs_{7}=0$. Using \cite[p. 81]{Toda}, we derive $[\e_{13},\e_{13}]=P(\e_{27})=\s\theta$ and $\ca_{10}\cg_{11}=\lt_{10}+\ch_{10}$. 

Therefore, $A$ is generated by the following two elements,  
$$(\ca_{13}\vee \cg_{13})\hc \mathrm{j}_{1}\hc \w_{14}\nu_{30}=\ca_{13}\w_{14}\nu_{30}=(\s\theta)\ch_{25},$$ 
$$(\ca_{13}\vee \cg_{13})\hc[\mathrm{j}_{1},\mathrm{j}_{2}]\hc\e_{33}=[\ca_{13},\cg_{13}]=[\e_{13},\e_{13}]\hc\ca_{25}\cg_{26}=(\s\theta)(\lt_{25}+\ch_{25}).$$
That is, $A=\mathrm{span}\n(\s\theta)\ch_{25},\,(\s\theta)(\lt_{25}+\ch_{25})\nn=\mathrm{span}\n(\s\theta)\ch_{25}\,, (\s\theta)\lt_{25}\nn.$ 
Here, we freely used some properties of Whitehead products, see \cite[Theorem  8.18,\,p. 484]{GW}).   However, substituting $n=0$ into the formula in Proposition\,\ref{ind}, we derive $\mathrm{Ind}=\ca_{13}\hc\pi_{33}(S^{14})+\cg_{13}\hc\pi_{33}(S^{20})=\mathrm{span}\n \ca_{13}\w_{14}\nu_{30}\nn=\mathrm{span}\n(\s\theta)\ch_{25}\nn.$
Through \cite[p. 52]{20STEM}, we have $$\pi_{33}(S^{13})=\mathrm{span}\n \ct_{13},\,(\s\theta)\ch_{25},\,(\s\theta)\lt_{25},\,\s\beta''' \nn\tg\z/8\jia(\z/2)^{\jia3}.$$ 
So $\mathrm{Ind}$ is not isomorphic to $A$, which tells us that for $n=0$ 
the formula  in Proposition\,\ref{ind} is not true in general.\end{proof}

\begin{rem}
    The matrix Toda bracket defined in this paper is a special case of the 3-fold one. All fundamental properties concerning 3-fold Toda brackets in \cite[p. 10-12]{Toda}, also hold for the matrix Toda bracket under certain conditions.
Specifically, for our matrix Toda bracket to hold valid, the indexing number $n$ is positive  and each wedge factor of the domain of  the formal matrices  engaged in  addition ``+"
is a double suspension; furthermore, when involving a formal  matrix taking ``the left composition with a homotopy class $\ck$ " (i.e., 
 $A\hc \ck$, where A is a formal matrix),  $\ck$ is a suspension. These assertions can be straightforwardly confirmed through direct verification.

\end{rem}

Whenever the Toda brackets are well-defined,
we also denote the Toda bracket
$\n a_{1}, \s^{n}(b_{1,1},b_{1,2}), \s^{n}(c_{1},c_{2})^{T} \nn_{n}$
as $\mtodaa{a_{1}}{\s^{n}  b_{1,1}     }{\s^{n}  c    _{1}}
{\s^{n}  b_{2,1}        }{\s^{n}  c    _{2}}_n$, called a right matrix Toda bracket; similarly, we denote $\n (a_{1},a_{2}), \s^{n}(b_{1,1},b_{2,1})^{T}, \s^{n}  c _{1}\nn$ as $\mtoda{  a  _{1}}{\s^{n}  b_{1,1}      }{\s^{n}  c    _{1}}{  a  _{2}}{\s^{n}  b_{2,1}     }_{n},$ called a left matrix Toda bracket.

There is a special case  of the Toda brackets in Definition \ref{lwndy}. That is, the  brackets of the form $$\left\{ (a_{1}, a_{2}), 
\;\,\s^{n}\begin{pmatrix}
b_{1,1}&   b_{1,2}        \\
b_{2,1} & 0      \\
  
\end{pmatrix},\;\, \s^{n}\begin{pmatrix}
c_{1}\\
c_{2} \end{pmatrix}
\right\}_{n},$$ ($\ell=r=2,\; b_{2,2}=0$), which in the case $n=0$ correspond to  \textbf{box brackets} up to sign  defined by Hardie, Marcum and  Oda in \cite{Oda2001}. Such  Toda brackets possess numerous advantageous properties.

 The Toda brackets $\left\{ (a_{1}, a_{2}), 
\;\,\s^{n}\begin{pmatrix}
b_{1,1}&   b_{1,2}        \\
b_{2,1} & 0      \\
  
\end{pmatrix},\;\, \s^{n}\begin{pmatrix}
c_{1}\\
c_{2} \end{pmatrix}
\right\}_{n}$ following our definition  can  be regarded as a practical model of   box brackets given by \cite{Oda2001} for explicit computations of unstable homotopy groups.
 Therefore, it is not unexpected that several properties of  box brackets  appear in our setting as well, see \cite[Proposition 3.2]{Oda2001} and the lemma below.
 
\begin{lem}  In Definition \ref{lwndy}, taking $\ell=r=2,\; b_{2,2}=0$,  further suppose  $X_{1}, X_{2}$ are suspensions and $n\geq1$. Then the  four conclusions stand valid.
    \begin{itemize}
        \item[\rm(1)] Substituting $a_{2}=0$, the matrix Toda bracket becomes a right matrix Toda bracket, that is,
        $$\left\{ (a_{1},  0), 
\;\,\s^{n}\begin{pmatrix}
b_{1,1}&   b_{1,2}        \\
b_{2,1} & 0      \\
  
\end{pmatrix},\;\, \s^{n}\begin{pmatrix}
c_{1}\\
c_{2} \end{pmatrix}
\right\}_{n}= 
\mtodaa{  a  _{1}}{\s^{n}  b_{1,1}     }{\s^{n}  c    _{1}}
{\s^{n}  b_{1,2}        }{\s^{n}  c    _{2}}_n.
$$
        \item[\rm(2)]Substituting $c_{2}=0$, the matrix Toda bracket becomes a left matrix Toda bracket, that is,
        $$\left\{ (a_{1}, a_{2}), 
\;\,\s^{n}\begin{pmatrix}
b_{1,1}&   b_{1,2}        \\
b_{2,1} & 0      \\
  
\end{pmatrix},\;\, \s^{n}\begin{pmatrix}
c_{1}\\
0 \end{pmatrix}
\right\}_{n}=\mtoda{  a  _{1}}{\s^{n}  b_{1,1}      }{\s^{n}  c    _{1}}{  a  _{2}}{\s^{n}  b_{2,1}     }_{n}.$$
        \item[\rm(3)] Substituting $b_{2,1}=0$ yields the following equation,
        $$\left\{ (a_{1}, a_{2}), 
\;\,\s^{n}\begin{pmatrix}
b_{1,1}&   b_{1,2}        \\
0 & 0      \\
  
\end{pmatrix},\;\, \s^{n}\begin{pmatrix}
c_{1}\\
c_{2} \end{pmatrix}
\right\}_{n}=  a  _{2}\hc\s^{n}[\s  U, X_{2}]+\mtodaa{  a  _{1}}{\s^{n}  b_{1,1}     \;\,}{\s^{n}  c    _{1}}
{\s^{n}  b_{1,2}       \;}{\s^{n}  c    _{2}}_n.$$
        \item[\rm(4)]  Substituting $b_{1,2}=0$ yields the following equation,
        $$\left\{ (a_{1}, a_{2}), 
\;\,\s^{n}\begin{pmatrix}
b_{1,1}&   0        \\
b_{2,1} & 0      \\
  
\end{pmatrix},\;\, \s^{n}\begin{pmatrix}
c_{1}\\
c_{2} \end{pmatrix}
\right\}_{n}=[\s^{n+1}Y_{2}, W]\hc\s^{n+1}  c    _{2}
+\mtoda{  a  _{1}}{\s^{n}  b_{1,1}      }{\s^{n}  c    _{2}}{  a  _{2}}{\s^{n}  b_{2,1}     }_{n}.$$
    \end{itemize}
\end{lem}

\begin{proof}  Referencing Proposition \ref{ind}, we delineate the formula employed to compute the indeterminacies. Let the maps $q_{k}: X_{1}\vee X_{2}\rightarrow X_{k}$,\,($k=1,2$) be the pinch maps.
\begin{itemize}
    \item[\rm(1)]By Definition \ref{lwndy} and  assumptions, we have\begin{eqnarray}
D&:= &\notag \left\{ (a_{1},  0), 
\;\,\s^{n}\begin{pmatrix}
b_{1,1}&   b_{1,2}        \\
b_{2,1} & 0      \\
  
\end{pmatrix},\;\, \s^{n}\begin{pmatrix}
c_{1}\\
c_{2} \end{pmatrix}
\right\}_{n}\\\notag
&= &\notag \left\{ (a_{1},  0), 
\;\,\s^{n}\begin{pmatrix}
b_{1,1}&   b_{1,2}        \\
0 & 0      \\
  
\end{pmatrix},\;\, \s^{n}\begin{pmatrix}
c_{1}\\
c_{2} \end{pmatrix}
\right\}_{n}+\left\{ (a_{1},  0), 
\;\,\s^{n}\begin{pmatrix}
0&   0       \\
b_{2,1} & 0      \\
  
\end{pmatrix},\;\, \s^{n}\begin{pmatrix}
c_{1}\\
c_{2} \end{pmatrix}
\right\}_{n}\\\notag
&= &\notag \left\{ a_{1} \hc\s^{n}q_{1},
\;\,\s^{n}\begin{pmatrix}
b_{1,1}&   b_{1,2}        \\
0 & 0      \\
  
\end{pmatrix},\;\, \s^{n}\begin{pmatrix}
c_{1}\\
c_{2} \end{pmatrix}
\right\}_{n}+\left\{ a_{1} \hc \s^{n}q_{1}, 
\;\,\s^{n}\begin{pmatrix}
0&   0       \\
b_{2,1} & 0      \\
  
\end{pmatrix},\;\, \s^{n}\begin{pmatrix}
c_{1}\\
c_{2} \end{pmatrix}
\right\}_{n}\\\notag
&\xyd &\notag  \left\{ a_{1}, 
\;\,\s^{n}(
b_{1,1},  b_{1,2}),\;\s^{n}\begin{pmatrix}
c_{1}\\
c_{2} \end{pmatrix}
\right\}_{n}+\left\{ a_{1}, 
\;\,\s^{n}(0,0),\;\s^{n}\begin{pmatrix}
c_{1}\\
c_{2} \end{pmatrix}
\right\}_{n}:=F'\\\notag
&\dyd &\notag \left\{ a_{1}, 
\;\,\s^{n}(
b_{1,1},  b_{1,2}),\;\s^{n}\begin{pmatrix}
c_{1}\\
c_{2} \end{pmatrix}
\right\}_{n}:=F.
\end{eqnarray}
\noindent   
Straightforward calculations confirm that  $\mathrm{Ind}( D )=\mathrm{Ind}( F' )=\mathrm{Ind}( F )$. So, $ D = F.$

     \item[\rm(2)]   By Definition \ref{lwndy} and  assumptions, we have\begin{eqnarray}
Q&:= &\notag \left\{ (a_{1}, a_{2}), 
\;\,\s^{n}\begin{pmatrix}
b_{1,1}&   b_{1,2}        \\
b_{2,1} & 0      \\
  
\end{pmatrix},\;\, \s^{n}\begin{pmatrix}
c_{1}\\
0 \end{pmatrix}
\right\}_{n}\\\notag
&= &\notag \left\{ (a_{1}, a_{2}), 
\;\,\s^{n}\begin{pmatrix}
b_{1,1}&   0        \\
b_{2,1} & 0      \\
  
\end{pmatrix},\;\, \s^{n}( \,\mathrm{j}_{1}  c_{1}\,)
\right\}_{n}+\left\{ (a_{1}, a_{2}), 
\;\,\s^{n}\begin{pmatrix}
0&   b_{1,2}        \\
0 & 0      \\
  
\end{pmatrix},\;\, \s^{n}( \,\mathrm{j}_{1}  c_{1}\,)
\right\}_{n}\\\notag
&\xyd &\notag \left\{ (a_{1}, a_{2}), 
\;\,\s^{n}(\begin{pmatrix}
b_{1,1}&   0        \\
b_{2,1} & 0      \\
  
\end{pmatrix}\hc  \mathrm{j}_{1} ),\;\, \s^{n}  c_{1} 
\right\}_{n}+\left\{ (a_{1}, a_{2}), 
\;\,\s^{n}(\begin{pmatrix}
0&   b_{1,2}        \\
0 & 0      \\
  
\end{pmatrix}\hc \mathrm{j}_{1} ),\;\, \s^{n}  c_{1} 
\right\}_{n}\\\notag
&= &\notag  \left\{ (a_{1}, a_{2}), 
\;\,\s^{n}\begin{pmatrix}
b_{1,1}       \\
b_{2,1}       \\
  
\end{pmatrix},\;\, \s^{n}  c_{1} 
\right\}_{n}+\left\{ (a_{1}, a_{2}), 
\;\,\s^{n}\begin{pmatrix}
0       \\
0       \\
  
\end{pmatrix},\;\, \s^{n}  c_{1} 
\right\}_{n}:=V'\\\notag
&\dyd  &\notag  \left\{ (a_{1}, a_{2}), 
\;\,\s^{n}\begin{pmatrix}
b_{1,1}       \\
b_{2,1}       \\
  
\end{pmatrix},\;\, \s^{n}  c_{1} 
\right\}_{n}:=V.
\end{eqnarray}

\noindent By direct verification, we have $\mathrm{Ind}(Q)=\mathrm{Ind}(V')=\mathrm{Ind}(V)$.
Thus, $Q=V$.

      \item[\rm(3)] By Definition \ref{lwndy} and  assumptions, we have\begin{eqnarray}
L&:= &\notag \left\{ (a_{1}, a_{2}), 
\;\,\s^{n}\begin{pmatrix}
b_{1,1}&   b_{1,2}        \\
0 & 0      \\
  
\end{pmatrix},\;\, \s^{n}\begin{pmatrix}
c_{1}\\
c_{2} \end{pmatrix}
\right\}_{n}\\\notag
&= &\notag \left\{ a_{1}\hc\s^{n}q_{1}, 
\;\,\s^{n}\begin{pmatrix}
b_{1,1}&   b_{1,2}        \\
0 & 0      \\
  
\end{pmatrix},\;\, \s^{n}\begin{pmatrix}
c_{1}\\
c_{2} \end{pmatrix}
\right\}_{n}+\left\{ a_{2}\hc\s^{n}q_{2}, 
\;\,\s^{n}\begin{pmatrix}
b_{1,1}&   b_{1,2}        \\
0 & 0      \\
  
\end{pmatrix},\;\, \s^{n}\begin{pmatrix}
c_{1}\\
c_{2} \end{pmatrix}
\right\}_{n}\\\notag
&\xyd &\notag \left\{ a_{1}, 
\;\,\s^{n}(
b_{1,1},  b_{1,2}),\;\, \s^{n}\begin{pmatrix}
c_{1}\\
c_{2} \end{pmatrix}
\right\}_{n}+\left\{ a_{2}, 
\;\,\s^{n}(
0,0),\;\, \s^{n}\begin{pmatrix}
c_{1}\\
c_{2} \end{pmatrix}
\right\}_{n}:=K\\\notag
&\dyd &\notag \left\{ a_{1}, 
\;\,\s^{n}(
b_{1,1},  b_{1,2}),\;\, \s^{n}\begin{pmatrix}
c_{1}\\
c_{2} \end{pmatrix}
\right\}_{n}:=M.
\end{eqnarray}
\noindent It is easy to infer $\ind (L)=\ind(K) $ which suggests $L=K$.
 Taking into account  an abelian  group $(G,+)$ with a subgroup $N$, if a set $S$ satisfies $G\dyd S\in G/N$, then  $S=S'+N$ for each subset $S'\xyd S$. Obviously, $\mathrm{Ind}(K)=\mathrm{Ind}(M)+  a  _{2}\hc\s^{n}[\s  U, X_{2}]$.
So, $K$ is equal to  
\begin{eqnarray}
\notag 
M+\mathrm{Ind}(K)&=& M+(\mathrm{Ind}(M)+  a  _{1}\hc\s^{n}[\s  U, X_{1}]) \\\notag
&=&(M+\mathrm{Ind}(M)))+  a  _{2}\hc\s^{n}[\s  U, X_{2}] \\\notag 
&=& M+  a  _{1}\hc\s^{n}[\s  U, X_{1}] \\\notag
&=&  a  _{1}\hc\s^{n}[\s  U, X_{1}]+M.
\end{eqnarray}
\item[\rm(4)] By synthesizing the methodologies of (2) and (3), it is easy to infer (4). The detail is left to the reader.
\end{itemize}
\end{proof}

\begin{rem}
 In this paper,  we only study the unstable  May-Lawrence matrix Toda bracket  which can be defined by the \textit{3}-fold Toda bracket,  does not investigate the unstable  May-Lawrence matrix Toda bracket  defined by the \textit{n}-fold, ($n\geq4$). 
 
 In fact, in the unstable range, even the \textit{4}-fold Toda bracket exhibits extreme complexity with various manifestations, (see \cite[Definition 4.12, p. 34-35]{dd}). Furthermore, both \cite{Ogt} and \cite{M}, the original references defining the \textit{4}-fold Toda brackets, contain mathematical errors on the unstable  \textit{4}-fold. These errors have been corrected by \cite{dd}, see \cite[p. 35-36]{dd} for details.
\end{rem}

  \indent With the current techniques, in  general we cannot ascertain the relationship between   the matrix Toda bracket indexed by 0 and the box bracket defined in $\textit{Top}_{*}$ given by  \cite{Oda2001}. What is clear is that in the stable range, our matrix Toda bracket of the form $\left\{ (a_{1}, a_{2}), 
\;\,\s^{n}\begin{pmatrix}
b_{1,1}&   b_{1,2}        \\
b_{2,1} & 0      \\
  
\end{pmatrix},\;\, \s^{n}\begin{pmatrix}
c_{1}\\
c_{2} \end{pmatrix}
\right\}_{n}$ and the box bracket in \cite{Oda2001} have same behaviours up to sign. Both left matrix Toda brackets following the definition in \cite[Section 4]{32STEM}  and this paper satisfy all properties given by \cite[Lem.\,4.1]{32STEM} up to sign if the indexing number $n\geq1$ and all the spaces $X_{k}$ are suspensions. Of course, they have  same behaviors up to sign in the stable range.

\section{The relative Homotopy and the \text{EHP} sequence} 
In the following definition, we will use the relative homoptopic relation $``\sim"$ of maps between pairs, (it is occasionally denoted by ``$\simeq$'' in some literatures). This can be found from some standard materials, for example,  \cite[p. 15] {Hsz} and      \cite[p. 8]{GW}.  We point out that in this sense, for two maps $f:(X,A)\rightarrow(X',A')$ and  $g:(X,A)\rightarrow(X',A')$ satisfying  $f\sim g \;\mathrm{rel} \;A$,  the condition $f|A=g|A$    is not required. The sets  of relative homotopy classes are studied by Toda and Mimura in \cite[Section 5]{Sp} and \cite{Toda1}.
\begin{definition}\label{tgxd} Let $A$ be a subspace of $X$, and $Z$ be a space.  The  set of relative homotopy classes  $[\s Z, (X,A) ]$ is defined as $[(C Z,Z), (X,A) ]$, that is, $$[\s Z, (X,A) ]:=\n f\, |\;\, f: (CZ,Z) \rightarrow(X,A) \nn/\sim,$$ where 
    $``\sim"$ denotes the relative homoptopic relation from  $(C Z,Z) $ to $ (X,A) $.
\end{definition}

It is well-known that $S^{1}$ is a \textit{co-H} space and its classical comultiplication $\mu'_{S^{1}}$  is the composition $$S^{1}\twoheadrightarrow S^{1}/S^{0}\stackrel{\cong}\longrightarrow S^{1}\vee S^{1}.$$ And $\s Z=Z\wedge S^{1}$ is  a \textit{co-H} space with the  classical comultiplication $\mu'=\mu'_{\s Z}= \mathrm{id}_{Z}\wedge\mu'_{S^{1}}$. Notice that $\mu'$  extends to the reduced cone. That is, we have the commutative diagram, where $C^{\mu'}=\mu'\wedge \mathrm{id}_{I}$ is   extended over $\mu'$, \\\indent \qquad \qquad \qquad \qquad \qquad\qquad \qquad \qquad 
\xymatrix{
  \s Z \ar[d]_{\xyd} \ar[r]^{\mu'} & \s Z\vee \s Z \ar[d]^{\xyd} \\
  C\s Z \ar[r]^{C^{\mu'}\quad} & C\s Z\vee C\s Z.  }
  
\noindent Similar to the absolute case,  
$$C^{\mu'}: (C\s Z , \s Z)\rightarrow(C\s Z\vee C\s Z, \s Z\vee \s Z) $$ gives an addition $``\,+\,$" on $[\s^{2}Z, (X,A) ]$  and gives a group structure on it. \\\indent The following lemma is essentially  from \cite[p. 422,\,429,\,430]{Toda1}, since there is the obvious natural  homeomorphism ($I$ with basepoint 1): $$C\s Z= (\s Z)\wedge I= (Z\wedge S^{1})\wedge I\cong (Z\wedge I)\wedge S^{1}=\s CZ.$$ 

\begin{lem} \label{xl}
    \begin{itemize}
        \item [\rm(1)] $[\s^{2}Z, (X,A)]$ is a group.
        \item [\rm(2)] $[\s^{3}Z, (X,A)]$ is an abelian group.
        \item [\rm(3)] Let $p:(C\s^{n}Z,\s^{n}Z)\rightarrow(\s^{n+1}Z,*)$  be the pinch map.  Then,  
         $$p^{*}:[(C\s^{n}Z,\s^{n}Z), (X,*)]\longrightarrow[(\s^{n+1}Z,*), (X,*)]=[\s^{n+1}Z, X]$$  are group isomorphisms for all $n\geq1$. Following from Definition \ref{tgxd}, we  denote above isomorphisms conveniently:\\\centerline{$[\s^{n+1}Z, (X,*)]\underset{\tg\;}{\stackrel{p^{*}}\longrightarrow}[\s^{n+1}Z, X],\;(n\geq1);$}
         in the case $n=0$,  $[\s^{}Z, (X,*)]\stackrel{p^{*}}\longrightarrow[\s^{}Z, X]$ is an isomorphism of pointed sets.
        \item [\rm(4)]Let  $i:A\hookrightarrow X$ and  $j:(X,*)\hookrightarrow (X,A)$ be the inclusions. Identify $[\s^{n+1}Z, (X,*)]$ with $[\s^{n+1}Z, X],\;(n\geq0)$. Then there exists a long exact sequence of groups \\\indent\qquad\qquad\qquad\qquad\xymatrix@C=0.34cm{
  \cdots\ar[r] &[\s^{2}Z, A] \ar[r]^{i_{*}} & [\s^{2}Z, X] \ar[r]^{j_{*}\quad} &[\s^{2}Z, (X,A)]  \ar[r]^{\quad\pa} & [\s Z, A] \ar[r]^{i_{*}} &[\s Z, X] ,  }

 and  an exact sequence of pointed sets\\ \indent \qquad\qquad\qquad\qquad\qquad
  \xymatrix@C=0.34cm{
  [\s Z, X] \ar[r]^{j_{*}\quad} &[\s Z, (X,A)]  \ar[r]^{\quad\pa} & [Z, A] \ar[r]^{i_{*}} &[Z, X] ,  }
  
 where
the boundary momorphism $\pa$ is given by $\pa[f]=[\,f|_{\s^{n}Z}\,]$\; for $f:(C\s^{n}Z,\s^{n}Z)\rightarrow (X,A)$, ($n\geq0$).
    \end{itemize}
\end{lem}
$\bzd$\\\\
\indent
Now we can show Theorem  \ref{hao1}. 

\noindent \textbf{Proof of Theorem 1.}
     Let spaces be localized at 2. We do something to change the homotopy fibration by  a strict fibration.  Let  \vspace{-1\baselineskip}$$N=JS^{m}\times_{H_{2}}(JS^{2m})^{I}=\n (x, f)\,| \;H_{2}(x)=f(0),\; x\in JS^{m},\; f: I\rightarrow JS^{2m} \;\text{is a free path} \nn  $$ \text{and}  \vspace{-1\baselineskip}$$\quad F=JS^{m}\times_{H_{2}}\mathrm{Path}(JS^{2m})=\n (x, f)\,| \;H_{2}(x)=f(0),\; x\in JS^{m},\; f: I\rightarrow JS^{2m}, f(1)=* \nn.$$
    Here, $I$ has basepoint 1,  and the inclusion $JS^{m}\stackrel{i}\hookrightarrow N$ is given by $x\mapsto (x,\;f_{_{H_{2}(x)}})$ where $f_{_{H_{2}(x)}}$ is the constant path whose image is $\n H_{2}(x) \nn$.\;Therefore, the homotopy fibration $JS^{m}\stackrel{H_{2}}\longrightarrow JS^{2m}\;\;(\text{with homotopy  fibre }S^{m})$ is changed by the strict  fibration $N\stackrel{\overline{H_{2}}}\longrightarrow JS^{2m}\;\; (\text{with   fibre}\; F),$  where $\overline{H_{2}}$  is the extension of $H_{2}$,  ($\overline{H_{2}}\hc i=H_{2}$). Thus, there exists a homotopy equivalence $S^{m}\stackrel{e}\longrightarrow F$. We know $H_{2}(x)=*$ for any $x\in S^{m}=J_{1}S^{m}\xyd JS^{m}$. Then we  obtain an inclusion $$(JS^{m}, S^{m})\stackrel{i}\hookrightarrow (N, F).$$  By checking the \textit{mod\;2} homology of dimension $m$, we know  the restriction $[\,i\,|_{S^{m}} ] $ is a generator of the $\z_{(2)}$-module $\pi_{m}(F)\tg\z_{(2)}$. Meanwhile,
      $[e]\in\pi_{m}(F)$ is also a generator, so $[\,i\,|_{S^{m}} ]=x[e]=[e]\hc x\e_{m}$ where  $x\in \z_{(2)}$ is not divisible by 2. Since $x(\mathrm{id}_{S^{m}})$ and $e$ are both homtopy equivalences, so    $i\,|_{S^{m}} : S^{m}\rightarrow F$ is also a homtopy equivalence. Then, by \textit{the Five Lemma} (note: \textit{the Five Lemma} still holds for the category of groups, but doesn't hold for the category of pointed sets), we have \vspace{-0.8\baselineskip}$$[\s^{r}Z,(JS^{m}, S^{m})]\stackrel{i_{*}}\longrightarrow [\s^{r}Z,(N, F)]$$
are isomorphisms of groups for all $r\geq2$.  Since $N\stackrel{\overline{H_{2}}}\longrightarrow JS^{2m}$ (with   fibre $F$) is a strict fibration, then $$\overline{H_{2}}_{*}:[\s^{r}Z, (N,F)]\rightarrow[\s^{r}Z, (JS^{2m},*)]$$ are isomorphisms of groups for all  $r\geq2$, (see \cite[Corollary\,8.8,\;p.  187]{GW} and \cite[Section 5]{Sp}). Since  $\overline{H_{2}}_{*}\hc i_{*}=H_{2*}$, then we derive the desired result.$\bzd$

\;\\ \indent For the above proposition, it is obvious that $H_{2*}$ can further become isomorphisms of $\z_{(2)}$-modules.
We  only need to notice for any $k\in\z_{(2)}$ and any  $\af :(C\s^{r-1}Z, \s^{r-1}Z ) \rightarrow (JS^{m}, S^{m})$, ($r\geq2$),  $$H_{2*}(k\af)=H_{2}\hc \af\hc k(\mathrm{id}_{C\s^{r-1}Z})=H_{2*}( \af)\hc k(\mathrm{id}_{C\s^{r-1}Z})=kH_{2*}(\af)$$
where  $k(\mathrm{id}_{C\s^{r-1}Z})=C^{k(\mathrm{id}_{\s^{r-1}Z})}=C^{(\mathrm{id}_{\s^{r-2}Z})\wedge k\e_{1}}$.

In the following proposition, we provide a clear and explicit construction of the boundary homomorphism $P$ in the generalized \textit{EHP} sequence and demonstrate its naturality.
\begin{pro}\label{gnlehp}
     After localization at 2, let the  sequence (called the generalized \textit{EHP} sequence) \\\\\indent\qquad\qquad\qquad\qquad\xymatrix@C=0.4cm{
 \cdots \ar[r]^{H\;\quad}& [\s^{2} Z,S^{m}]\ar[r]^{\s\;\;\,\;} & [\s^{3}Z,S^{m+1}] \ar[r]^{H\;} &[\s^{3}Z,S^{2m+1}]\ar[r]^{\;\;\;P} & [\s Z,S^{m}] \ar[r]^{\s\;\;\,\;} & [\s^{2}Z,S^{m+1}]   }
\;\\\\ be given  in which each $H$ is induced by $H_{2}$ up to the classical isomorphism, and each $P\colon[\s^{r+1}Z, S^{2m+1}]\rightarrow [\s^{r-1}Z, S^{m}]$ is  defined to be $P=\pa\hc (H_{2*})^{-1}\hc \lo_{1}$. Here $\lo_{1}: [\s\,-\;, S^{2m+1}]\stackrel{\is}\longrightarrow [\,-,\;JS^{2m}]$ is the classical isomorphism and  $\pa: [\s^{r}Z, (JS^{m}, S^{m})]\rightarrow [\s^{r-1}Z,  S^{m}]$  is the  boundary homomorphism  defined  by our Lemma   \ref{xl}\,(3), ($r\geq2$).  Then,  this sequence is exact in the category of groups.  Moreover, with respect to the category of  2-local suspensions and maps, $$\mbox{\Large\{}
[\s^{2}Z',S^{2m+1}]\stackrel{P}\rightarrow [ Z',S^{m}] \;{\Large | }\; Z'\; \text{is a 2-local suspension}\,\mbox{\Large\}}$$
is a natural transformation from the contravariant founctor $[-,\;S^{m}]$ to the contravariant founctor  $[\s^{2}-,\;S^{2m+1}]$.

\end{pro}
    
    \begin{proof}
        Let spaces be localized at 2. By Theorem \ref{hao1} and Lemma   \ref{xl}\,(4), we know that the generalized \textit{EHP} sequence is exact in the category of groups. Identify $\lo_{1}$ with $\lo_{0}:[\s\,-\,,\;S^{2m+1}]\stackrel{\is}\longrightarrow[\,-\,,\;\lo S^{2m+1}]$. Recall from \cite[(1.12),\,p. 8]{Toda} that $\beta\hc \s\af=\lo_{0}^{-1}((\lo_{0}\beta)\hc\af)$. And notice that  the boundary homomorphism $\pa$ is natural by its definition. Since $P=\pa\hc (H_{2*})^{-1}\hc \lo_{1}$, then $P$ is natural.
\end{proof}

\begin{rem}As in common knowledge,   we can derive a long exact sequence by applying the functor $[\s Z, -]$ to the 2-local homotopy fibration sequence, $\cdots\longrightarrow\lo^{2} S^{2m+1}\rightarrow S^{m}\hookrightarrow \lo S^{m+1}\stackrel{ H_{2}}\longrightarrow\lo S^{2m+1}$. Proposition \ref{gnlehp} primarily highlights the $\mathbf{specific\; constructions}$ of the boundary homomorphisms $P$, in order to use \cite[Proposition 2.4]{I} to derive the $H$-formula.
\end{rem}

\indent By using of the naturality of  homomorphisms $P$ of the  generalized \textit{EHP} sequence (Proposition\,\ref{gnlehp}), we have the following corollary.

\begin{cor}\label{zrx}

After localization at 2, if $L_{1},L_{2}$ are suspensions, then
     there exists a  commutative diagram   where $p_{1}:L_{1}\vee L_{2}\rightarrow L_{1}$ is the pinch map and  $i_{2}: L_{2}\rightarrow L_{1}\vee L_{2}$ is the inclusion map. Moreover,  the vertical arrows are short exact sequences which both  split into direct products and the split homomorphisms are induced by the  pinch and inclusion maps. \\ \\\indent \qquad\qquad\qquad\qquad\qquad\qquad\qquad
\xymatrix{
  [\s^{2}L_{1}, S^{2m+1}]\ar[d]_{(\s^{2}p_{1})^{*}} \ar[r]^{\qquad P} &[L_{1}, S^{m}] \ar[d]^{p_{1}^{*}} \\
 [\s^{2}(L_{1}\vee L_{2}), S^{2m+1}]\ar[d]_{(\s^{2}i_{2})^{*}}\ar[r]^{\qquad P} &[L_{1}\vee L_{2}, S^{m}]\ar[d]^{i_{2}^{*}} \\
  [\s^{2}L_{2}, S^{2m+1}]\ar[r]^{\qquad P} &[ L_{2}, S^{m}].   }

 \end{cor}
 $\bzd$\\
 
 \indent Next, we give  the proof of Theorem \ref{zdl2}.
 
\noindent \textbf{Proof of Theorem 2.} By Proposition \ref{gnlehp} and Corollary \ref{zrx},  we infer $$P^{-1}(\mathbbm{f}_{1\times\ell}\hc \mathbbm{b}_{\ell\times r})=\sum_{s=1}^{r}(\s^{2}p_{s})^{*}P^{-1}(\mathbbm{f}_{1\times\ell}\hc\mathbbm{b}_{\bullet,s})$$ where $p_{s}: \bigvee_{j=1}^{r} Y_{j}\rightarrow Y_{s} $  ($1\leq s\leq r$) are the pinch maps. 
    Then the theorem follows from Definition \ref{lwndy} and
    \cite[Proposition 2.4]{I}. Notice that in  \cite[Proposition 2.4]{I},\; $h=H_{2}$ and
 $\Gamma\hc\s^{2}=P^{-1}\hc\s^{2}$ where the homomorphism $P$ follows the definition given by Proposition \ref{gnlehp}.$\bzd$

\section{An Application}
In this section, all spaces, maps and homotopy classes    are localized at 2. The homotopy group $\pi_{i}(S^{n})$ localized at 2 is also denoted by $\pi_{i}^{n}.$ \\\indent 
 To show the following proposition, we make some  preparations here,\\ \vspace{-0.7\baselineskip}\\\centerline{$ \cg'\nu\ty\nu_{7}\cg_{10}\m2\nu_{7}\cg_{10},\;$ (\cite[(7.19),\,p. 71]{Toda})$,\;\ca\ch=\nu^{3}$,} \\ \centerline{$\nu_{5}\cg_{8}\nu_{15}^{2}=\ca_{5}\bar{\lt}_{6}$, (\cite[Lemma   12.10, p. 143]{Toda}) and  $\bar{\lt}_{6}=\ca_{6}\cn_{7}$, (\cite[(10.23),\,p. 110]{Toda}).}
\\\noindent Thus, \vspace{-2\baselineskip}\\\;\begin{equation}
\cg'\ca_{14}\ch_{15}=\cg'\nu_{14}\hc\nu_{17}^{2}=\nu_{7}\cg_{10}\nu_{17}^{2}=\ca_{7}\bar{\lt}_{8}=\ca_{7}^{2}\cn_{8}\tag{5.1} \label{ybgs}    
\end{equation}

\begin{pro}\label{51}

In $\pi_{26}^{6}$, the generator $\ct'$ is taken as $$
\ct'\in\mtodaa{\nu_6}{\eta_9}{\ca_{10}\cn_{11}}
{\s^{2}(\cg'\ca_{14})\;\;}{\bar{\nu}_{17}}_2
.$$
Then, $H(\ct')=\ca_{11}\cn_{12}.$  All properties of the original $\ct'$ in \cite{20STEM}  satisfies, this new choice of $\ct'$ still satisfies.   
\end{pro}

\begin{proof}
We denote the matrix Toda bracket shown in our proposition  by $T$.
We list the following facts from  \cite[p. 44, p. 50, p. 70]{Toda}:$$\nu_{5}\ca_{8}=P(\e_{11}),\;\s^{2}\cg'=2\cg_{9},\;P^{-1}(\nu_{5} \hc \s(\cg'\ca_{14}))\xyd \pi_{18}^{11}=\mathrm{span}\n\cg_{11}\nn.$$  
Then, $\nu_{6}\ca_{9}=\nu_{6}\hc \s^{2}(\cg'\ca_{14})=0.$   Formula (\ref{ybgs}) gives
$\ca_{7}\hc\ca_{8}\cn_{9}+\cg'\ca_{14}\hc\bar{\nu}_{15}=0$. Hence, $T$ is well-defined. We know $\cg_{11}\ch_{18}=0,$ (\cite[(10.18),\,p. 108]{Toda}).
Then, by Theorm\;\ref{zdl2}, 
$$H(T)\xyd -(P^{-1}(\nu_{5}\ca_{8}))\hc\ca_{11}\cn_{12}-(P^{-1}(\nu_{5} \hc \s(\cg'\ca_{14})))\hc\ch_{18}
=\ca_{11}\cn_{12}.$$
Since $\pi_{18}^{6}=\mathrm{span}\n[\e_{6},\e_{6}]\hc \cg_{11}\nn$, (\cite[p. 74]{Toda}),
then,  by Definition\,\ref{lwndy}, 
\begin{eqnarray}
\notag T& 
\xyd& \n\nu_{6}, \;\ca_{9}\vee0_{9}^{(8)}, \;\mathrm{j}_{1}\ca_{10}\cn_{11} +\mathrm{j}_{2}\ch_{17}\nn_{0}  \\\notag&
\xyd& \n\nu_{6}, \;\ca_{9}\vee0_{9}^{(8)}, \;\mathrm{j}_{1}\ca_{10}\cn_{11} \nn+\n\nu_{6}, \;\ca_{9}\vee0_{9}^{(8)}, \;\mathrm{j}_{2}\ch_{17}\nn   \\\notag&
\xyd&\n\nu_{6}, \;\ca_{9}, \;\ca_{10}\cn_{11} \nn+\n\nu_{6}, 0_{9}^{(8)}, \;\ch_{17}\nn  \\\notag&
=& \n\nu_{6}, \;\ca_{9}, \;\ca_{10}\cn_{11} \nn+\pi_{18}^{6}\hc \ch_{18}   \\\notag&
=& \n\nu_{6}, \;\ca_{9}, \;\ca_{10}\cn_{11} \nn. 
\end{eqnarray}
(Here, $0_{9}^{(8)}=0\in\pi_{9+8}^{9}$).\;By the proof of \cite[Lemma   15.3,\,p. 43]{20STEM}, we see that any element in $\n\nu_{6}, \;\ca_{9}, \;\ca_{10}\cn_{11} \nn$ can be chosen as $\ct'$. Thus, we  take $\ct'\in T$. And so, $H(\ct')=\ca_{11}\cn_{12}$. 
\end{proof}\;\\
\indent We see  \cite[p. 47]{20STEM} only suggests $H(\ct')\ty \ca_{11}\cn_{12}\m8\s^{2}\co'$.\;Hence our proposition gives a more sophisticated construction of $\ct'$. \\\indent Notice $H\pi_{6+5}^{6}=\mathrm{span}\n H(P(\e_{13}))\nn=\mathrm{span}\n 2\e_{11}\nn$, (\cite[Proposition\,2.7, Proposition\,5.9]{Toda}). Following from \cite[Proposition\,2.2]{Toda}, $H(\pi_{6+5}^{6}\hc\ca_{11}\cn_{12})=0$. We know $\pi_{9+17}^{9}=\s\pi_{8+17}^{8}$,\;(\cite[(12.9),\,p. 146]{Toda}). Then
$\n\nu_{6}, \;\ca_{9}, \;\ca_{10}\cn_{11} \nn$ is a coset   of $$\pi_{6+5}^{6}\hc\ca_{11}\cn_{12}+\nu_{6}\hc\pi_{9+17}^{9}\xyd \s\pi_{5+20}^{5}.$$ Recall from the proof of Proposition\,\ref{51} that $$T\xyd \n\nu_{6}, \;\ca_{9}, \;\ca_{10}\cn_{11} \nn,\; H(T)=\ca_{11}\cn_{12};$$ and notice $H\hc \s=0$. These lead to the following result.
\begin{cor}
    $H\n\nu_{6}, \;\ca_{9}, \;\ca_{10}\cn_{11} \nn= \ca_{11}\cn_{12}.$$\bzd$
\end{cor}

\end{document}